%% file: BFECC_uniform_grid_manuscript.tex
\documentclass[review,hidelinks,onefignum,onetabnum]{siamart220329}

\input{ex_shared}

\ifpdf
\hypersetup{
  pdftitle={A ghost-point based second order accurate finite difference method on uniform orthogonal grids for electromagnetic scattering around PEC},
  pdfauthor={H. Lee and Y. Liu}
}
\fi

\usepackage{bm, float}

\usepackage{lipsum}
\usepackage{array}
\usepackage{subfig} 
\usepackage{amsfonts}
\usepackage{graphicx}
\usepackage{epstopdf}
\usepackage{algorithmic}
\ifpdf
\DeclareGraphicsExtensions{.eps,.pdf,.png,.jpg}
\else
\DeclareGraphicsExtensions{.eps}
\fi

\begin{document}

\maketitle

\begin{abstract}
	We propose a finite difference method to solve Maxwell's equations in time domain in the presence of a perfect electric conductor that impedes the propagations of electromagnetic waves. Our method is a modification of the existing approach by Zou and Liu \cite{zl21}, from a locally perturbed body-fitted grid to a uniform orthogonal one for complicated PEC objects. Similar to their work we extrapolate ghost point values by exploiting the level set function of the interface and the PDE-based extension technique,  which allows us to circumvent scrutinizing local geometries of the interface. We stipulate a mild requirement on the accuracy of our extrapolation that the ghost values need only be locally second order accurate. Nevertheless the resulting accuracy of our method is second order thanks to the application of back and forth error correction and compensation, which also relaxes CFL conditions. We demonstrate the effectiveness of our approach with some numerical examples. 
\end{abstract}



\section{Introduction}
Electromagnetic scattering problems have long been a popular subject of numerical simulations due to their wide range of applications including remote sensing, antenna and optical designs, to name a few. As is often encountered in the realm of computational physics, numerical treatment of material interfaces is also a challenging task when it comes to solving Maxwell's equations in the presence of scattering objects across which material properties may exhibit a sharp transition. One such prototypical example is concerned with a perfect electric conductor (PEC), an idealized material with infinite conductivity. The resulting PEC boundary conditions necessarily involve the boundary geometries which may not always be simple. Moreover the normal and tangential components of electric and magnetic fields, respectively, are discontinuous across the boundary, while they are both identically zero within the PEC. The aforementioned two issues, along with practical considerations, make it imperative to devise an accurate and efficient numerical method to capture the PEC boundary conditions. In this work we propose a simple, yet second order accurate finite difference method for the scattering of electromagnetic waves around the PECs. Our approach benefits largely from the level set framework \cite{o88} and the back and forth error compensation and correction (BFECC) \cite{dupontBackForthError2003}. The former has garnered much attention across a wide span of scientific communities, and among a multitude of research work we refer to \cite{o01} for a review. The latter was originally designed to enhance the interface level set and fluid simulations \cite{k05,SFKLR08}, but has recently been expanded in its applications to computational electrodynamics \cite{wangBackForthError2019,zl21}. We also refer to the modified MacCormack method \cite{SFKLR08} which is unconditionally stable.

Before we turn to our method we must clearly delineate that the finite-difference time domain method (FDTD), more commonly known as Yee scheme \cite{yeeNumericalSolutionInital1966}, has long been a mainstay in electrodynamic simulations since its first conception over 50 years ago. The popularity of Yee scheme can succinctly be attributed to its ``simplicity, efficiency, and generality'' \cite{r99}; we refer to \cite{k93} for the details of FDTD. Arguably the simplest FDTD approach to boundary treatment is a staircased approximation \cite{umashankarCalculationExperimentalValidation1987a} but its accuracy can deteriorate to the first order \cite{DDH01}. As alternatives globally or locally nonorthogonal conformal grids \cite{jin-fa-leeModelingThreedimensionalDiscontinuities1992,fangLocallyConformedFinitedifference1993} have been used to more accurately enforce boundary conditions while there has also been a line of work on overlapping grids  \cite{y92,liuOverlappingYeeFDTD2009}. However such modifications of finite difference grids in one form or another are inextricably tied to the issue of mesh generation, which can become very costly especially in high dimensions. The authors in \cite{DDH01} adhere to the conventional orthogonal grid, yet locally modifies the finite difference stencil to achieve second order accuracy. In \cite{ln21} a fourth order accurate finite difference method on orthogonal grids is proposed based on the correction function method which entails a minimization problem. After all we focus in this work on closely resembling the aforementioned merits of FDTD while preserving its second order accuracy. 

The main contribution of our work is to propose a numerical method that extends the previous work \cite{zl21} on irregular grids to uniform rectangular grids. The authors therein rely on locally conforming point-shifted grids for the interface which is implicitly represented by some level set function. Then the PDE extension technique \cite{AS98,Peng99,Fedkiw99} is applied to facilitate construction of ghost points, which is done automatically without explicitly taking into account the boundary geometry. Once the ghost points are established, they employ the first order numerical scheme followed by the application of BFECC. The latter provides the accuracy improvement from the first order to the second order in addition to increased CFL conditions. 

We seek uniform rectangular grids so that our proposed numerical method can readily be applied to various PEC geometries in different problems. On one hand, the simplicity of the uniform grids is highly desirable in practice and comes with ease of implementation. Moreover there is no need to perform the point-shift algorithm which could run into problems when the interface
is locally underresolved as noted by the authors in \cite{zl21}. 
On the other hand, our reasoning is that construction of ghost points should in principle be independent of whether they match with the interface; one can ultimately interpolate the ghost values based on case-by-case analysis of how the boundary geometry intersect the grid points. Provided that level set functions are available there is no need for additional information about the interfaces. The consequence of adopting uniform grids is that the grid points adjacent to the interfaces now carry necessary boundary data when the PDE extension techniques are harnessed. Once the ghost points are set up, however, our method relies on the underlying first order numerical schemes and BFECC, hence fully inheriting the simplicity and robustness of the existing method on the point-shifted grid. What distinguishes our method is that it is supported by theoretical underpinnings of the BFECC which have thus far been rigorously proved only on uniform grids.

The rest of the paper is organized as follows. We present our proposed method in detail in Section \ref{sec:method}. We then demonstrate performance of our method by means of numerical experiments in Section \ref{sec:num_experiments}. A conclusion is given in Section \ref{sec:conc}.

\section{Method}
\label{sec:method}

We discuss the various ingredients of our proposed method, highlighting the similarities and differences between ours and the existing algorithm on the shifted grid \cite{zl21}.

\subsection{Back and Forth Error Compensation and Correction}
\label{Sec:BFECC}

Similar to the shifted grid case we incorporate the Back and Forth Error Compensation and Correction method, or BFECC \cite{dupontBackForthError2003}. To recapitulate the idea let us consider a linear constant coefficients first order hyperbolic PDE system of the form:
\begin{align*}\label{eq:hyperbolic_system_forward}
\partial_t \bm{u} + \sum_{i = 1}^{d} A_i \partial_{x_i} \bm{u} = 0,
\end{align*}
and its forward in time numerical approximation
\[
\bm{U}^{n+1} = \mathcal{L} \bm{U}^{n}.
\]
Here $\bm{U}^{n+1}$ denotes the numerical solution at time $t_{n+1}$ obtained by applying a linear numerical scheme to advance $\bm{U}^{n}$ to the next time step. Then we march the solution $\bm{U}^{n+1}$ backward in time to obtain
\[
\widetilde{\bm{U}}^{n} = \mathcal{L}^{\star} \bm{U}^{n+1}
\]
where $ \mathcal{L}^{\star}$ represents the resulting numerical scheme obtained by applying the same scheme $\mathcal{L}$ to the time-reversed equation of  (\ref{eq:hyperbolic_system_forward}). Now we define the back and forth error as 
\[
\bm{e}^n =  \frac{1}{2} \left( \bm{U}^n - \tilde{\bm{U}}^n \right),
\]
which accounts for numerical errors committed in one forward time stepping. It is by this error $\bm{e}^{n}$ that we compensate the solution $\bm{U}^n$ to correct $\bm{U}^{n+1}$. That is, we recompute $\bm{U}^{n+1}$ by
\[
\bm{U}^{n+1} = \mathcal{L} \left( \bm{U}^n + \bm{e}^{n} \right).
\]
It has been proved \cite{dupontliu07} that on uniform grids BFECC improves the $r$-th order accurate linear numerical scheme to $(r+1)$-th order, where $r$ is odd. Hence we only need to design a first order  scheme to guarantee second order accuracy. In terms of stability BFECC is proved to stabilize numerical schemes that would otherwise be unstable for all CFL numbers such as the forward in time and centered in space numerical scheme \cite{dupontliu07}. Moreover the $L^2$ stability is attainable with CFL numbers that are larger than 1 since BFECC permits the Fourier growth factor of no more than $2$ for the underlying scheme; see Theorem 1 in \cite{dupontliu07} for scalar equations and Theorem 3.1 in \cite{wangBackForthError2019} for systems.

Let us next present our numerical scheme $\mathcal{L}$ tailored
to the present study of the following two dimensional transverse magnetic Maxwell's equations
\begin{align}\label{eqn:2d_maxwell_tmz}
\begin{split}
& \frac{\partial H_x}{\partial t} = - \frac{\partial E_z}{\partial y} \\
& \frac{\partial H_y}{\partial t} =  \frac{\partial E_z}{\partial x} \\
& \frac{\partial E_z}{\partial t} = \frac{\partial H_y}{\partial x} - \frac{\partial H_x}{\partial y}
\end{split}
\end{align}
where both the permittivity $\epsilon$ and permeability $\mu$ parameters are normalized to 1. The equations are to be complemented with appropriate far-field boundary conditions as well as the PEC boundary conditions which we will discuss in sequel. In their work   \cite{wangBackForthError2019} the authors numerically simulate Maxwell's equations using a parameterized $\theta$-scheme $\mathcal{L}_\theta:= (1-\theta) \mathcal{C} + \theta \mathcal{F}$ as the underlying scheme for BFECC, where $\mathcal{C}$ and $\mathcal{F}$ denotes the centered and Lax-Friedirchs schemes, respectively. We choose $\theta = \frac{4}{5}$ as done in \cite{wangBackForthError2019,zl21}, thereby obtaining in the following discretizations for  (\ref{eqn:2d_maxwell_tmz}):
\begin{enumerate}
	\item{\bf Forward}. \\
	\begin{eqnarray} \label{eq:bf}
	&(\widetilde{E_z})^{n+1}_{i,j} = \displaystyle \frac{\left( ({E_z})^{n}_{i-1,j} + ({E_z})^{n}_{i,j} + ({E_z})^{n}_{i+1,j} + ({E_z})^{n}_{i,j-1} + ({E_z})^{n}_{i,j+1}\right)}{5} +  \nonumber \\
	& \displaystyle \frac{\Delta t}{2 \Delta x} \left((H_y)^{n}_{i+1,j} - (H_y)^{n}_{i-1,j}\right) - \displaystyle \frac{\Delta t}{2 \Delta y} \left((H_x)^{n}_{i,j+1} - (H_x)^{n}_{i,j-1}\right),
	\end{eqnarray} and analogously for equations of $H_x$ and $H_y$.
	\item{\bf Backward}. \\
	\begin{eqnarray}  \label{eq:bfb}
	&(\widetilde{E_z})^{n}_{i,j} = \displaystyle \frac{\left( (\widetilde{E_z})^{n+1}_{i-1,j} + (\widetilde{E_z})^{n+1}_{i,j} + (\widetilde{E_z})^{n+1}_{i+1,j} + (\widetilde{E_z})^{n+1}_{i,j-1} + (\widetilde{E_z})^{n+1}_{i,j+1}\right)}{5} - \nonumber \\
	& \displaystyle \frac{\Delta t}{2 \Delta x} \left((\widetilde{H_y})^{n+1}_{i+1,j} - (\widetilde{H_y})^{n+1}_{i-1,j}\right) + \displaystyle \frac{\Delta t}{2 \Delta y} \left((\widetilde{H_x})^{n+1}_{i,j+1} - (\widetilde{H_x})^{n+1}_{i,j-1}\right),
	\end{eqnarray} and analogously for equations of $H_x$ and $H_y$.
	\item{\bf Error compensation}. Reassign the value
	\begin{equation} \label{eq:be}
	({E_z})^{n}_{i,j} = ({E_z})^{n}_{i,j} + \frac{1}{2}(({E_z})^{n}_{i,j} - (\widetilde{E_z})^{n}_{i,j})
	\end{equation}
	and analogously for $H_x$ and $H_y$.
	\item{\bf Final forward marching with the compensated solution at time $t_n$}. \\
	\begin{eqnarray} \label{eq:bfbff}
	&({E_z})^{n+1}_{i,j} = \displaystyle \frac{\left( ({E_z})^{n}_{i-1,j} + ({E_z})^{n}_{i,j} + ({E_z})^{n}_{i+1,j} + ({E_z})^{n}_{i,j-1} + ({E_z})^{n}_{i,j+1}\right)}{5} + \nonumber \\
	& \displaystyle \frac{\Delta t}{2 \Delta x} \left((H_y)^{n}_{i+1,j} - (H_y)^{n}_{i-1,j}\right) - \displaystyle \frac{\Delta t}{2 \Delta y} \left((H_x)^{n}_{i,j+1} - (H_x)^{n}_{i,j-1}\right),
	\end{eqnarray} and analogously for equations of $H_x$ and $H_y$.
\end{enumerate} 

The spatial accuracy of the underlying scheme \eqref{eq:bf} is first order due to the 5-point averaging. We introduce such averaging in order to hinder possible appearance of numerical artifacts around the corners of PECs. The analysis in \cite{wangBackForthError2019} establishes that the ($\theta$-dependent) CFL number lies between $\sqrt{1.5}$ and $\sqrt{2}$. This makes it possible to use larger time steps such as the one corresponding to CFL number $1$, which is prohibited in the classical FDTD framework due to numerical instabilities \cite{m01}.

\subsection{PEC boundary treatment via level set approach}
The PEC boundary conditions are given by
\begin{equation} \label{eq:pec}
\bm{E} \times \bm{n} = 0 \;\; {\rm and} \;\; \bm{H} \cdot \bm{n} = 0
\end{equation}
where $\bm{n}$ denotes the outward normal to the interface. In order to satisfy the boundary conditions we resort to the method of images \cite{jackson1999classical} by which ghost values are set up across the interface. More specifically we need to extrapolate the ghost values inside the PEC based on the field values outside the PEC so that the resulting fields are smooth and satisfy the boundary conditions (\ref{eq:pec}). Our uniform grid in general does not match the PEC geometries, hence the boundary conditions need to be enforced implicitly by those grid points in the vicinity of the interfaces.

To facilitate the implementation of the image method we follow \cite{zl21} wherein the level set framework is used to extend some field values $u$ ``off" the interface \cite{Peng99}. Let $\phi$ denote the signed distance function to the interface which takes positive and negative values inside and outside the interface, respectively. We will write with abuse of notation $\bm{n} = \frac{\nabla \phi}{|\nabla \phi|}$ everywhere in the domain. Then we apply the following transport equation 
\begin{align}\label{eq:level_set}
\begin{split}
u_t + \bm{n} \cdot \nabla u = 0~\\
\end{split}
\end{align} which we discretize as in the case of the Maxwell's equations to obtain:  
\begin{align}\label{eq:update_eqn_level_set}
\begin{split}
u^{n+1}_{i,j} &= \displaystyle \frac{\left( {u}^{n}_{i-1,j} + {u}^{n}_{i,j} + {u}^{n}_{i+1,j} + {u}^{n}_{i,j-1} + {u}^{n}_{i,j+1}\right)}{5} - \Delta t [\bm{n}]_{i,j} \cdot  [\nabla u]^{n}_{i,j},\\
\end{split}
\end{align}
where $[\nabla u]^{n}_{i,j}$ and $[\bm{n}]_{i,j}$ are the centered finite difference approximations of $\nabla u$ and $\nabla \phi$, respectively. Scheme (\ref{eq:update_eqn_level_set}) is iterated for sufficient number of times until the fields values in the first ghost layer are populated and reach steady states. We follow the work \cite{zl21} and set $\Delta t = 0.2 \Delta x$ in (\ref{eq:update_eqn_level_set}). In view of the boundary conditions (\ref{eq:pec}), the quantities $u$ that need to be extended are  $E_z$, $\bm{H} \cdot \bm{n}$ and $\bm{H} \cdot \bm{t}$, where $\bm{t} = R \bm{n}$ where $R$ is the $2 \times 2$ rotation matrix in the clockwise direction by $\frac{\pi}{2}$. Let us first discuss the signed distance function $\phi$ before we proceed further with details of our construction of fictitious fields values.

\subsubsection{Construction of signed distance function $\phi$}
An indispensable ingredient in our method is accurate computation of $\phi$ and $\nabla \phi$. The uniform grid used in our method does not capture any information about the interface so it is imperative to obtain an accurate implicit representation of the interface via $\phi$. While there are costs involved in computation of $\phi$, it is a one time expense only at the beginning of our numerical simulations. Moreover an accurate resolution of $\phi$ is required only near the interface in order to define the ghost values using the transport equation (\ref{eq:update_eqn_level_set}). As will be seen we need to approximate  $\phi$ and $\nabla \phi$ with at least second order accuracy in order to make sure that the ghost values are second order accurate, which in turn renders our underlying scheme such as (\ref{eq:bf}) and (\ref{eq:bfb}) first order accurate. Let us also note that the second order accuracy requirement on $\nabla \phi$  is attributed to the change of variables from $H_x$ and $H_y$ to $\bm{H} \cdot \bm{n}$ and $\bm{H} \cdot \bm{t}$ which should introduce no more than $O((\Delta x)^2)$ errors.

One way to obtain $\phi$ and $\nabla \phi$ with the desired accuracy is to consider re-distancing equation
until equilibrium \cite{sussmanLevelSetApproach1994}

\begin{align}\label{eqn:distance_level_set}
\begin{split}
\frac{\partial \phi}{\partial t} + \bm{v} \cdot \nabla \phi = sgn (\phi)~,\\
\end{split}
\end{align}
subject to some initial level set data. Here $\bm{v} = sgn(\phi) \frac{\nabla \phi}{| \nabla \phi |}$, and 
\begin{equation*}
sgn(x) =
\begin{cases}
1, \quad x > 0 \\
0, \quad x = 0 \\
-1, \quad x < 0 \\
\end{cases}.
\end{equation*} 
Equation (\ref{eqn:distance_level_set}) can be discretized following the work of \cite{m10} which is shown to yield third order accurate $\phi$ near smooth interface. We can then obtain by a simple centered difference second order accurate approximation of $\nabla \phi$ near the interface. As an alternative to the dynamic (\ref{eqn:distance_level_set}) one can consider the time independent Eikonal equation \cite{d13} and apply a myriad of numerical methods such as an improved version of the fast marching method \cite{ch01,zhao05}. Since $\phi$ needs to be computed once only at the initial time near the interface, we can afford the existing high order schemes or employ much more refined sub-grids to facilitate computation of $\phi$. Although our proposed method relies on accurate $\phi$, a key ingredient in our approach is the \emph{application} of $\phi$ in (\ref{eq:update_eqn_level_set}) to assign ghost values within the PEC. For clarity and simplicity of our presentation we curtail further detailed discussion of $\phi$, and assume from here on that a highly accurate $\phi$ is available at our disposal.

\subsubsection{Construction of odd images across the interface}
We need to enforce zero Dirichlet boundary conditions for $E_z$ and $\bm{H}\cdot\bm{n}$. To this end we follow the work of \cite{zl21} so that the smoothly extended field values across the interface are odd functions. Let us write $u(x_g)$ to denote the value of some field $u$ at a ghost point  location $x_g$ which is inside, yet $O(\Delta x)$ away from the PEC interface. Then a Taylor series expansion is given by
\[
u(x_g) = u(x_g^{\star}) + (x_g-x_g^{\star})\cdot \nabla u (x_g^{\star}) + O((\Delta x)^2)
\]
where $x_g^{\star}$ denotes the closest point on the interface to $x_g$. Then since $x_g-x_g^{\star} = \phi(x_g) \bm{n}(x_g)+O((\Delta x)^2)$ and $\phi(x_g)=O(\Delta x) $, imposing the homogeneous Dirichlet condition $u(x_g^{\star}) = 0$ yields
\[
u(x_g) = \phi(x_g) \bm{n}(x_g^{\star})\cdot \nabla u (x_g^{\star}) + O((\Delta x)^2).
\]
The key idea adopted in \cite{zl21} is then to extend the quantity $\bm{n} \cdot \nabla u$ ``off" the interface using the PDE extension (\ref{eq:update_eqn_level_set}). However our grid points in general do not fall exactly on the interface, so the information $\bm{n} \cdot \nabla u$ is not as readily available as in \cite{zl21}. To remedy this issue we propose that the missing information at the interface is first approximated just outside the PEC before it is extended across the interface to the first ghost layer over a distance of say $ 2\Delta x$ as opposed to $\Delta x$. That is, we proceed as follows:
\begin{enumerate}
	\item[(1)] Identify the set $S$ of grid points $(i,j)$ that are in the \emph{second} boundary layer. That is, $(i,j) \in S$ if $(i,j)$ as well as all its four neighbors $(i-1,j), (i+1,j), (i,j-1), (i,j+1)$ are outside the PEC, but at least one of its neighbors is adjacent to the PEC (a grid point is outside but adjacent to the PEC if at least one of its four neighbors
	is inside the PEC.) Then we compute for all $(i,j) \in S$
	$$
	w_{i,j} = \frac{u_{i,j}}{\phi_{i,j}}
	$$  which is a first order approximation of $\bm{n}\cdot\nabla u$ at the interface (since it's an odd extension.)
	\item[(2)] For all $(i,j)$ on $\bm{n}$ side of $S$ but not in $S$, apply scheme (\ref{eq:update_eqn_level_set}) (with $u$ replaced with $w$) to extend $w$ along $\bm{n}$ until the first boundary layer (grid points outside but adjacent to the PEC) and ghost layer are populated with the extended values, say  $\widetilde{w}_{i,j}$. Define the ghost values and redefine the first boundary layer by 
	\[
	u_{i,j} = \widetilde{w}_{i,j} \phi_{i,j}.
	\] We call the redefinition of the first boundary layer the linearization step in our algorithm.
\end{enumerate}
We need to clarify some aspects of the aforementioned steps in more details. In the first step we use the grid points in the second boundary layer since some grid points in the first boundary layer can be very close to the interface depending on the geometry of the PEC. This can in turn potentially lead to numerical instabilities due to very small $\phi_{i,j}$; even with the exact $\phi$ we observed instabilities in some of our numerical experiments when we extended $w_{i,j}$ across the interface when it was computed using the first boundary layer grid points. In fact a similar observation is also noted in some studies of particle methods as in \cite{m97} which proposes implementing a lower bound on the closest distance to the interface. In a similar spirit is the point shifted case \cite{zl21} where grid points within $\frac{\Delta x}{2}$ distance to the interface along a grid line are shifted to the interface, thus eliminating the situation. We note that the one-sided approximation ${u_{i,j}}/{\phi_{i,j}}$ is $O(\Delta x)$ accurate provided that  $\phi_{i,j}$ is second order accurate.

One can view the redefinition in the first boundary layer as a local second order approximation to exact solution values there because of the odd extension across the interface. The odd symmetry would be better enforced with respect to the constructed first ghost values and we observed in our numerical tests that this numerical technique provides better convergence rates without any numerical artifacts.

\subsubsection{Construction of even images across the interface}
So far as the tangential components $\bm{H} \cdot \bm{t}$ are concerned we need to extend our field values in an even manner so that zero Neumann boundary conditions $\frac{\partial (\bm{H} \cdot \bm{t})}{\partial \bm{n}} = 0$ are enforced. Our approach is to construct a quadratic polynomial approximation adjacent to the interface, similar in spirit to the work \cite{kpy04}, but based on the PDE extension (\ref{eq:level_set}). To illustrate the idea  assume for simplicity that the vector fields of normal vectors $\bm{n}$ are straight lines near the interface. We can then consider without loss of generality the one dimensional setting along a normal line with points $x_{i} < x_{k}$  where $x_{i}$ and $x_{k}$ belonging to the outside and inside of the interface, respectively. Let us denote by $x_{\star}$ the interface location between $x_{i}$ and $x_{k}$. Then a smooth extension of $u$ at the ghost point $x_{k}$ which satisfies the zero Neumann boundary condition is given by   
\[
\begin{aligned}
u(x_{k}) &= \frac{u^\prime(x_i)}{2(x_{i}-x_\star)}(x_{k}-x_\star)^2 + u(x_i) + \frac{u^\prime(x_i)}{2}(x_\star-x_i) \\
&= \frac{u^\prime(x_i)}{2\phi_{i}}(\phi_{k})^2 + u(x_i) - \frac{u^\prime(x_i)}{2}(\phi_i). 
\end{aligned}
\]
If we generalize this to our current two dimensional setting, then the two quantities to be constantly extended across the interface are the concavity term $\frac{\nabla u (x_{i,j})\cdot \bm{n}_{i,j}}{2\phi_{i,j}}$ and the constant term $u(x_{i,j})-\frac{\nabla u (x_{i,j})\cdot \bm{n}_{i,j}}{2} \phi_{i,j}$ (the extrapolated value at the PEC boundary.) We next turn to the case where the vector fields $\bm{n}$ are not straight lines adjacent to the interface. Then the errors in the normal directions $\bm{n}$ imply that the accuracy of our quadratic extension reduces to second order. However our underlying scheme still remains first order accurate, which hence permits the application of BFECC as in the odd extension case. 

In the case of the point-shifted grid \cite{zl21}, the authors use linear approximations for the even extensions of tangential components of the magnetic fields. We adopt here quadratic approximations which incurs no further computational complexities as our quadratic polynomials are effectively one dimensional along the normal directions $\bm{n}$. 
This also allows the situation in which the ghost grid points are very close to the PEC boundary. Note that this situation doesn't exist in the point-shifted grid used in \cite{zl21}.
For our even extension we proceed as follows.

\begin{enumerate}
	\item[(1)] For all $(i,j) $ in the second boundary layer $S$ defined as in the odd extension case, we compute
	$$
	v_{i,j} = \frac{[\bm{n}]_{i,j} \cdot  [\nabla u]^{n}_{i,j}}{2\phi_{i,j}} \quad \text{and} \quad 
	z_{i,j} = u_{i,j} - \frac{[\bm{n}]_{i,j} \cdot  [\nabla u]^{n}_{i,j}}{2} \phi_{i,j}.
	$$ 
	\item[(2)] For all $(i,j) \not\in S$, apply scheme (\ref{eq:update_eqn_level_set})(with $u$ replaced with $v$ or $z$) to extend $v$ and $z$ along local $\bm{n}$ until the first boundary and ghost layers are populated with the extended values, say  $\widetilde{v}_{i,j}$ and $\widetilde{w}_{i,j}$.   Define the values in the first ghost layer by 
	\[
	u_{i,j} = \widetilde{v}_{i,j} (\phi_{i,j})^2 + \widetilde{z}_{i,j}.
	\]
\end{enumerate} 
The second boundary layer used in the first step allows us to apply the centered difference schemes to compute the derivatives. As opposed to the odd extension case we choose not to redefine the values in the first boundary layer using our quadratic extension. This is because the redefined values necessarily involves approximation errors, and the solution is nonlinear across the interface.

In our numerical tests we indeed observed more accurate results when the values in the first boundary layer are not reconstructed using the quadratic extensions.

\subsection{Algorithm}
We can summarize our method in what follows. Note our construction of ghost values are carried out before each time stepping, both forward and backward. We also implement the post-processing (Step 3) to polish our numerical solution in the first boundary layer as it has been been treated as a ghost layer during the odd extensions throughout the entire duration of simulation

\renewcommand{\thealgorithm}{}
\begin{algorithm}[H]
	\caption{ }
	\begin{enumerate}
		\item Compute the signed distance $\phi$ with $\phi > 0$ in the interior of PEC. Compute the unit normal $\bm{n}$ and tangential $\bm{t}$ vectors using the centered differences.
		\item For $n = 0,1, \dots N-1$, where $T = N \Delta t$, do:
		\begin{enumerate}
			\item[(i)] Decompose $\bm{H}$ into $\bm{H} \cdot \bm{n}$ and $\bm{H} \cdot \bm{t}$.
			\item[(ii)] Extend $E_z,\bm{H} \cdot \bm{n},\bm{H} \cdot \bm{t}$ across the interface: the former two by odd extensions and the latter by an even extension.
			\item[(iii)] Reassemble $\bm{H} = (\bm{H} \cdot \bm{n}) \bm{n} + (\bm{H} \cdot \bm{t}) \bm{t}$. 
			\item[(iv)] Apply the forward advection step with scheme (\ref{eq:bf}).
			\item[(v)] Repeat the steps (i),(ii),(iii), and apply the backward advection step with scheme (\ref{eq:bfb}).
			\item[(vi)] Compute the back and forth error, and compensate the solution by (\ref{eq:be}).
			\item[(vii)] Repeat the steps (i),(ii),(iii), and apply the final forward advection step with scheme (\ref{eq:bfbff}).
		\end{enumerate}
		\item Do the steps (i),(ii),(iii) to redefine $E_z$ and $\bm{H}\cdot\bm{n}$, (hence $\bm{H}$), in the first boundary layer. Set the ghost values of all field variables to $0$.
	\end{enumerate}
	\addtocounter{algorithm}{-1}
\end{algorithm}

Before we turn to numerical experiments we should clarify what far field boundary conditions are incorporated in our method. We apply the unsplit convolutional Perfectly Matched Layers (PMLs)\cite{k07} on the outskirts of our computational domain as the absorbing boundary conditions. We closely follow the implementation in \cite{wangBackForthError2019} which demonstrates the interplay between the BFECC and the PMLs. As pointed out in \cite{a18}, analytical solutions can instead be used as outer boundary conditions when they are available. Alongside the PMLs we follow the total field/scattered field formulations \cite{ta05} as commonly done in many FDTD simulations.

\section{Numerical experiments}\label{sec:num_experiments}
\subsection{Perfect Electric Conductor with a circular cross-section}

We first consider the circular PEC centered at $(0.5,0.5)$ with the radius of $0.2$ in the computational domain $\Omega = [0,1]^2$.  We use the exact signed distance function $\phi$ and compute its gradients $\bm{n}$ using a centered finite difference scheme. We compute the solutions on the grid of size $\Delta x = \Delta y =\frac{1}{20}, \frac{1}{40}, \frac{1}{80}, \frac{1}{160}$, which are then compared with the reference solution computed on $\Delta x = \Delta y = \frac{1}{640}$. The numerical errors, measured in discrete $l_1$ norms, are calculated over the layer of the width $0.1$ outside the PEC boundary, $0<\phi<0.1$.

\subsubsection{A soliton Gaussian pulse}
Following the numerical example in \cite{ln21} we consider the Gaussian incident wave given by
\begin{align}
H_x(x,y,t) &= 0  \nonumber \\
H_y(x,y,t) &=  -\frac{(x-\gamma-t)}{\sigma^2}\exp\left(-\left( \frac{x-\gamma-t}{\sigma}\right)^2\right)\\
E_z(x,y,t) &= \frac{(x-\gamma-t)}{\sigma^2}\exp\left(-\left( \frac{x-\gamma-t}{\sigma}\right)^2\right) \nonumber
\end{align}
where $\sigma = 0.1, \gamma = -0.1$. We choose the terminal time of $T = 0.4$ so that the incident wave has reached the boundary of PEC. As can be seen in Table \ref{tb:e1gc} we obtain the second order rates of convergence.

\begin{table}[H]
	\centering
	\begin{tabular}{|c|c|c|c|c|c|c|c|c|c|c|} \hline
		& \multicolumn{2}{|c|}{ $E_z$ } & \multicolumn{2}{|c|}{  $B_x$ }& \multicolumn{2}{|c|}{ $B_y$ } \\ \hline
		$\Delta x = \Delta y$ & Error & Order & Error & Order & Error & Order \\ \hline
		$1/20$ & $   1.18 \times 10^{0} $ & -- & $4.18\times 10^{-1} 
		$ & -- & $7.23 \times 10^{-1}  $ & --\\ \hline
		$1/40$ & $     5.58 \times 10^{-1}$ &    $1.08$
		& $  2.03 \times 10^{-1}$ & $1.04$ & $ 3.92 \times 10^{-1}$ &     $0.88$ \\ \hline
		$1/80$ & $
		1.63 \times 10^{-1}$  & 
		$1.78$ & $   7.04 \times 10^{-2}$ & $1.53$ &    $1.51 \times 10^{-1}$  &     $1.37$
		\\ \hline
		$1/160$ & $
		3.88 \times 10^{-2}$  &
		$2.07$ & $   1.86 \times 10^{-2}$ & $1.92$ & $3.94 \times 10^{-2}$ &     $1.94$
		\\ \hline
	\end{tabular}
	\caption{Errors for the circular PEC with Gaussian pulse, $T = 0.4$, $\Delta t /\Delta x = 1$. }\label{tb:e1gc}
\end{table}

\begin{figure}[H]
	\centering
	\subfloat[Circular PEC]{\includegraphics[scale=0.4]{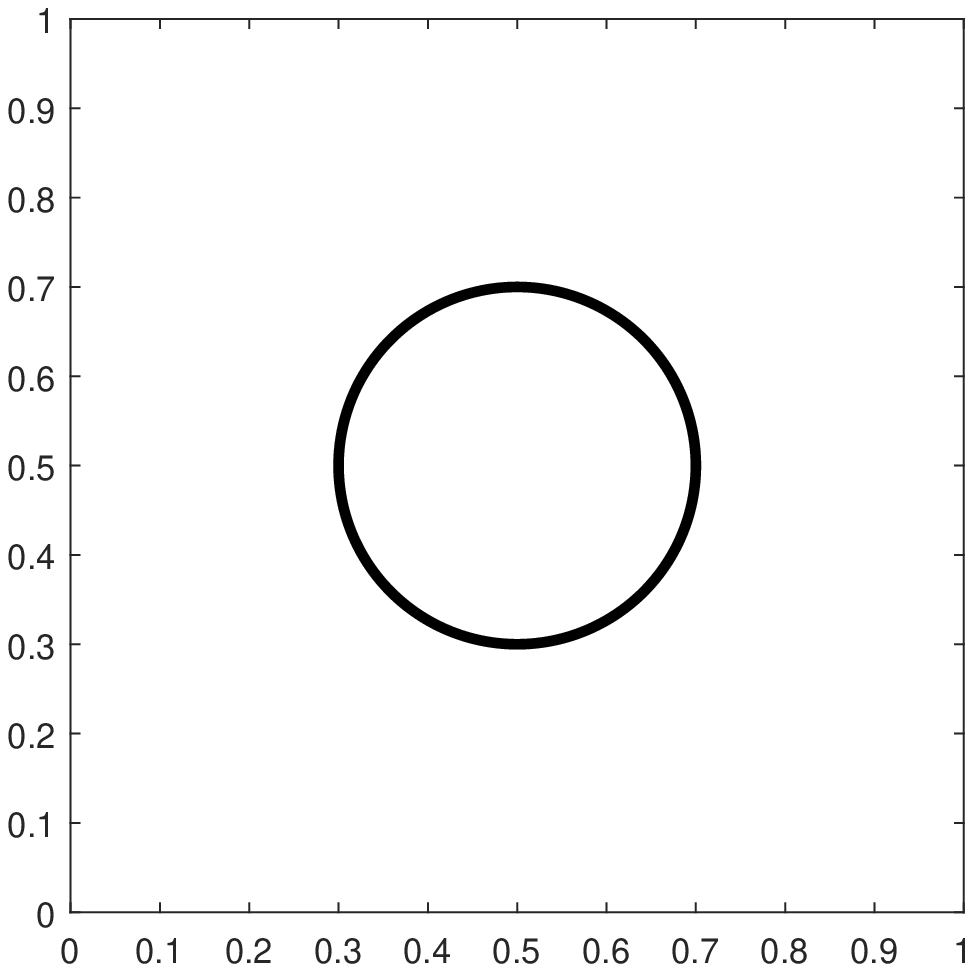}}
	\subfloat[$E_z$]{\includegraphics[scale=0.4]{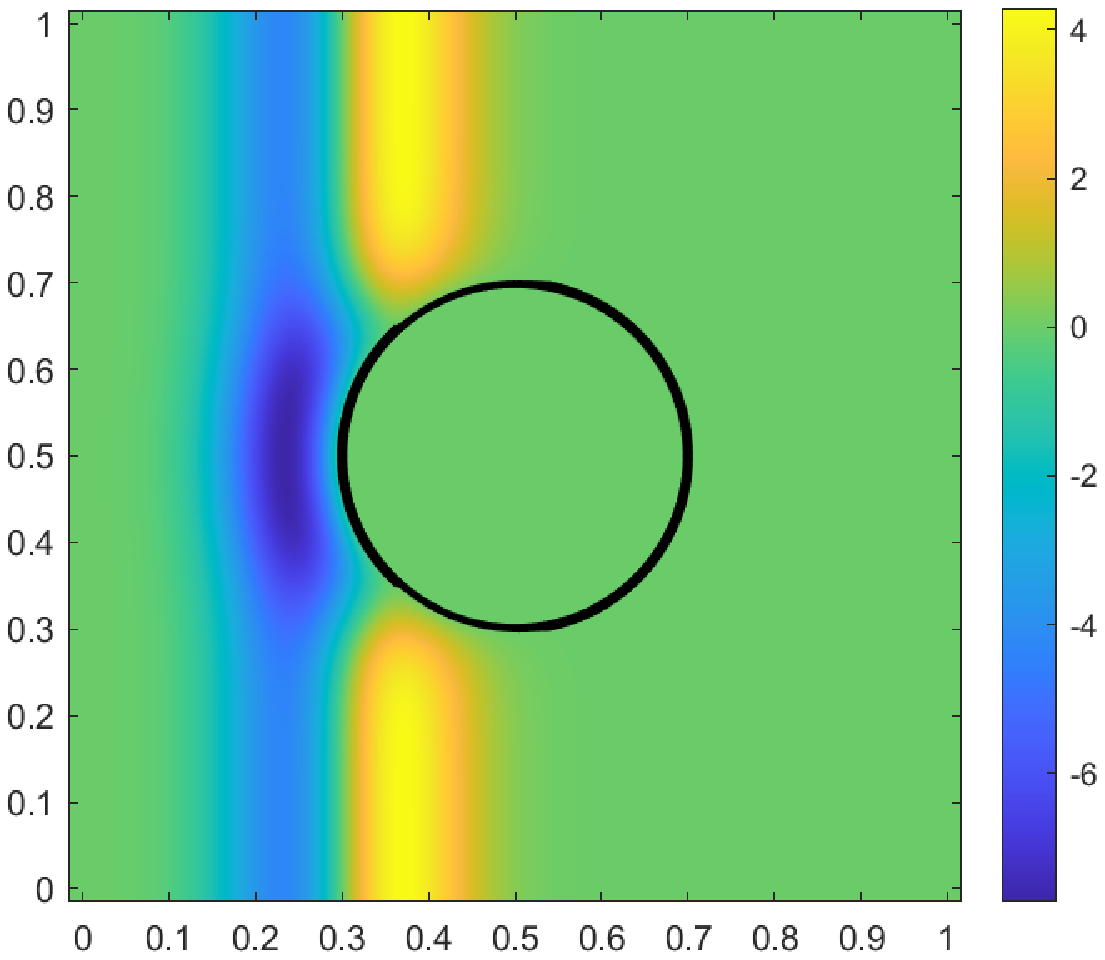}}
	\caption{Solution snapshot of $E_z$ at time $T = 0.4$ with the Gaussian incident wave.}  
\end{figure}

\subsubsection{A periodic plane wave}
Next we consider the incident wave similar to the one in \cite{m81}
\begin{align}
H_x(x,y,t) &= 0  \nonumber \\
H_y(x,y,t) &=  - \sin( \omega (x-t))) \chi(x-t)\\
E_z(x,y,t) &=  \sin( \omega (x-t)) \chi(x-t) \nonumber
\end{align} where  $\omega = 2\pi/0.3$ is the angular frequency and $\chi(z)$ denotes the heaviside step function. We introduce the step function in order to satisfy the zeroth order compatibility condition between the initial data and the boundary data \cite{f02}, in particular $ E_z(x,y,0) = 0$ on the boundary of the PEC. We compute the numerical solution until $T = 0.8$ so that the wavefront of the incident wave has been scattered by and diffracted around the PEC. We again observe the second order accuracy as in the single pulse case; see Table \ref{tb:e2gc}).

\begin{table}[H]
	\centering
	\begin{tabular}{|c|c|c|c|c|c|c|c|c|c|c|} \hline
		& \multicolumn{2}{|c|}{ $E_z$ } & \multicolumn{2}{|c|}{  $B_x$ }& \multicolumn{2}{|c|}{ $B_y$ } \\ \hline
		$\Delta x = \Delta y$ & Error & Order & Error & Order & Error & Order \\ \hline
		$1/20$ & $   5.13 \times 10^{-1} $ & -- & $3.25\times 10^{-1} 
		$ & -- & $4.31 \times 10^{-1}  $ & --\\ \hline
		$1/40$ & $     2.21 \times 10^{-1}$ &    $1.21$
		& $  1.58 \times 10^{-1}$ & $1.04$ & $ 2.16 \times 10^{-1}$ &     $1.00$ \\ \hline
		$1/80$ & $
		6.23 \times 10^{-2}$  & 
		$1.83$ & $   4.28 \times 10^{-2}$ & $1.88$ &    $5.46 \times 10^{-2}$  &     $1.98$
		\\ \hline
		$1/160$ & $
		1.47 \times 10^{-2}$  &
		$2.09$ & $   1.02 \times 10^{-2}$ & $2.07$ & $1.23 \times 10^{-2}$ &     $2.15$
		\\ \hline
	\end{tabular}
	\caption{Errors for the circular PEC with plane wave, $T= 0.8$, $\Delta t /\Delta x = 1$}\label{tb:e2gc}
\end{table}

\begin{figure}[H]
	\centering
	\includegraphics[scale=0.2]{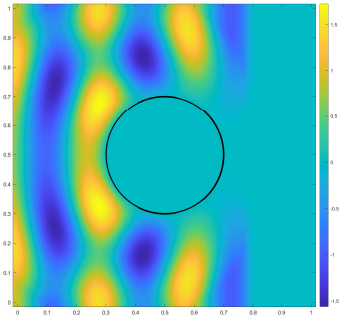}
	\caption{Solution snapshot of $E_z$ at time $T = 0.8$ with the sine incident wave.}
\end{figure}

\subsection{Cylindrical PEC with a re-entrant corner}
\begin{figure}[H]
	\centering
	\subfloat[Circular wedge PEC]{\includegraphics[scale=0.4]{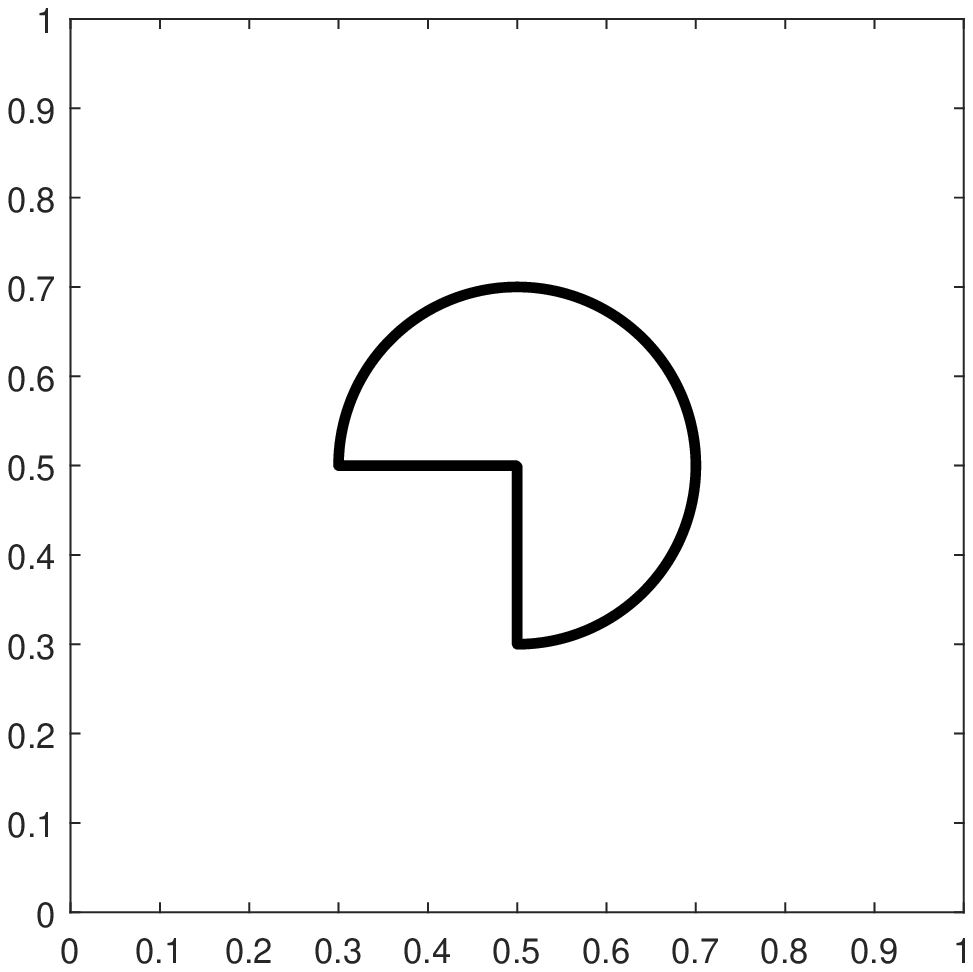}}
	\subfloat[$E_z$]{\includegraphics[scale=0.4]{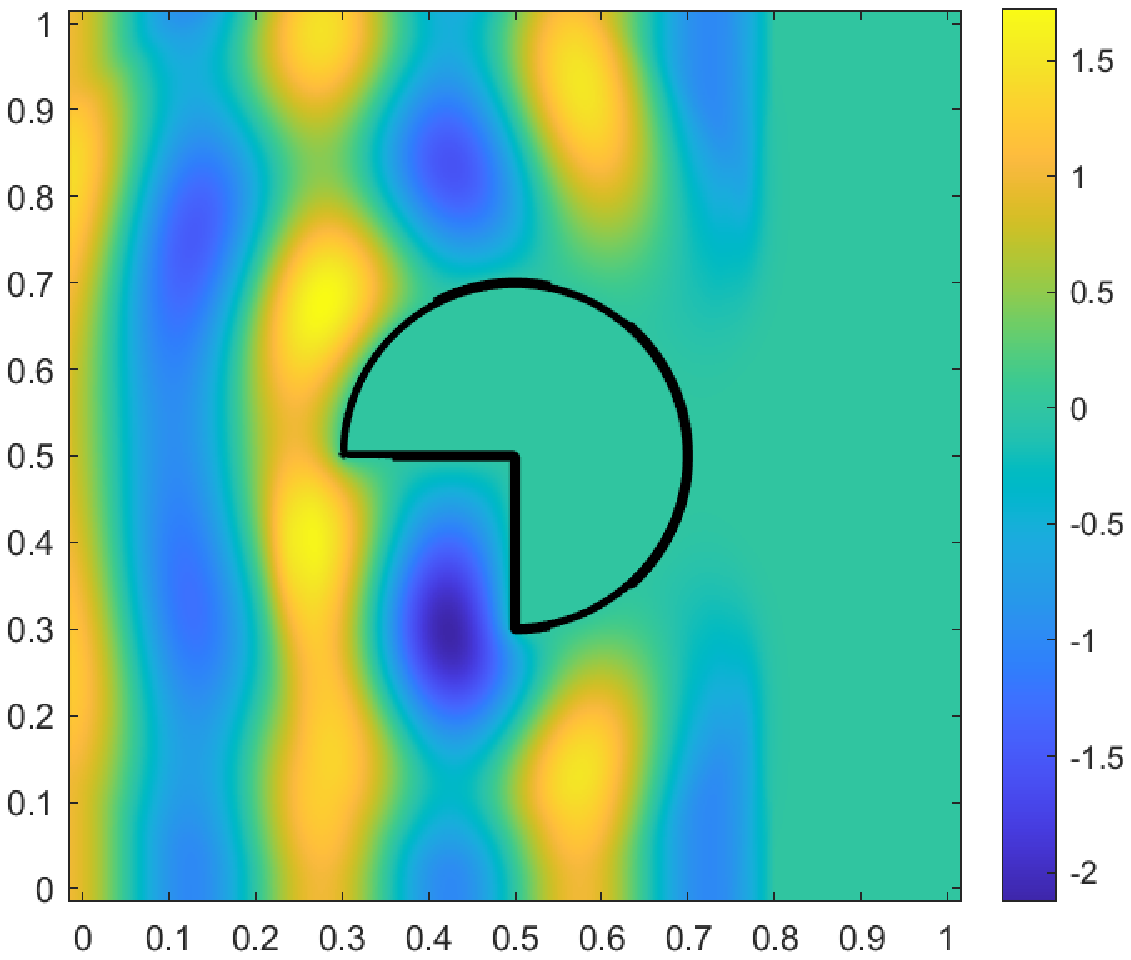}}
	\caption{Solution snapshot of $E_z$ at time $T = 0.8$ with the sine incident wave and circular wedge.}
	\label{fig:corner_PEC}
\end{figure}

We are interested in how our method performs when the PEC boundary has sharp corners. To illustrate we consider the circular arc cross section which subtends the angle of $\frac{3\pi}{2}$ as shown in Figure \ref{fig:corner_PEC}. The radii of the circular part of the PEC is set to $0.2$. We use the same incident plane wave as in the circular PEC case and compute the solutions till the terminal time $T = 0.8$. Despite the singularity of analytic solutions near the corner, the convergence rates stay close to the second order. Moreover we do not observe any numerical artifacts in our computed solutions.

\begin{table}[H]
	\centering
	\begin{tabular}{|c|c|c|c|c|c|c|c|c|c|c|} \hline
		& \multicolumn{2}{|c|}{ $E_z$ } & \multicolumn{2}{|c|}{  $B_x$ }& \multicolumn{2}{|c|}{ $B_y$ } \\ \hline
		$\Delta x = \Delta y$ & Error & Order & Error & Order & Error & Order \\ \hline
		$1/20$ & $   6.14 \times 10^{-1} $ & -- & $3.27\times 10^{-1} 
		$ & -- & $3.99 \times 10^{-1}  $ & --\\ \hline
		$1/40$ & $     2.37 \times 10^{-1}$ &    $1.37$
		& $  1.88 \times 10^{-1}$ & $0.79$ & $ 2.63 \times 10^{-1}$ &     $0.60$ \\ \hline
		$1/80$ & $
		6.65 \times 10^{-2}$  & 
		$1.83$ & $   5.43 \times 10^{-2}$ & $1.80$ &    $7.65 \times 10^{-2}$  &     $1.78$
		\\ \hline
		$1/160$ & $
		1.74 \times 10^{-2}$  &
		$1.93$ & $   1.46 \times 10^{-2}$ & $1.89$ & $1.99 \times 10^{-2}$ &     $1.94$
		\\ \hline
	\end{tabular}
	\caption{Errors for the $3/4$-moon shaped PEC with plane wave, $T= 0.8$, $\Delta t /\Delta x = 1$}\label{tb:e3gc}
\end{table}

\subsection{Two circular PECs with re-entrant corners}
{
	We next consider the case where there is more than one PEC object. It is not difficult to envision that the reflected electromagnetic waves end up interacting again by yet another PEC nearby. For concreteness we employ two identical copies of a circular PEC with the angle $3\pi/2$ as in Section 3.2. They are positioned at $(0.3,0.3)$ and $(0.6,0.6)$ with the radii of $0.15$. We use the same numerical set up as in the single corner case except with the modification of the angular frequency of the incident wave to $2\pi/0.2$. In addition to the corner singularities due to each individual PEC, the signed distance function has kinks at the equidistant points to the two PECs. Nevertheless our method is still shown to be quite effective as can be seen in Table 4. To better resolve the singularities, one may also consider introducing some nonlinear limiter as done in \cite{h16}, which we leave for future studies.

\begin{table}[H] 
	\centering
	\begin{tabular}{|c|c|c|c|c|c|c|c|c|c|c|} \hline
		& \multicolumn{2}{|c|}{ $E_z$ } & \multicolumn{2}{|c|}{  $B_x$ }& \multicolumn{2}{|c|}{ $B_y$ } \\ \hline
		$\Delta x = \Delta y$ & Error & Order & Error & Order & Error & Order \\ \hline
		$1/20$ & $   4.52 \times 10^{-1}   $ & -- & $ 2.58\times 10^{-1} 
		$ & -- & $3.85 \times 10^{-1}  $ & --\\ \hline
		$1/40$ & $     3.81\times 10^{-1}$ &    $0.25$
		& $  2.66 \times 10^{-1}$ & N.A & $ 4.07 \times 10^{-1}$ &     N.A \\ \hline
		$1/80$ & $
		1.78 \times 10^{-1}$  & 
		$1.10$ & $ 1.25 \times 10^{-2}$ & $1.09$ &    $1.76 \times 10^{-2}$  &     $1.21$
		\\ \hline
		$1/160$ & $
		5.86 \times 10^{-2}$  &
		$1.60$ & $  3.54  \times 10^{-2}$ & $1.82$ & $4.53 \times 10^{-2}$ &     $1.95$
		\\ \hline
	\end{tabular}
	\caption{Errors for the two $\frac{3}{4}$-moon shaped PEC with plane wave, $T= 0.8$, $\Delta t / \Delta x = 1$}
\end{table}
	
	\begin{figure}[H]
	\centering
	\subfloat[Two circular wedge PECs]{\includegraphics[scale=0.4]{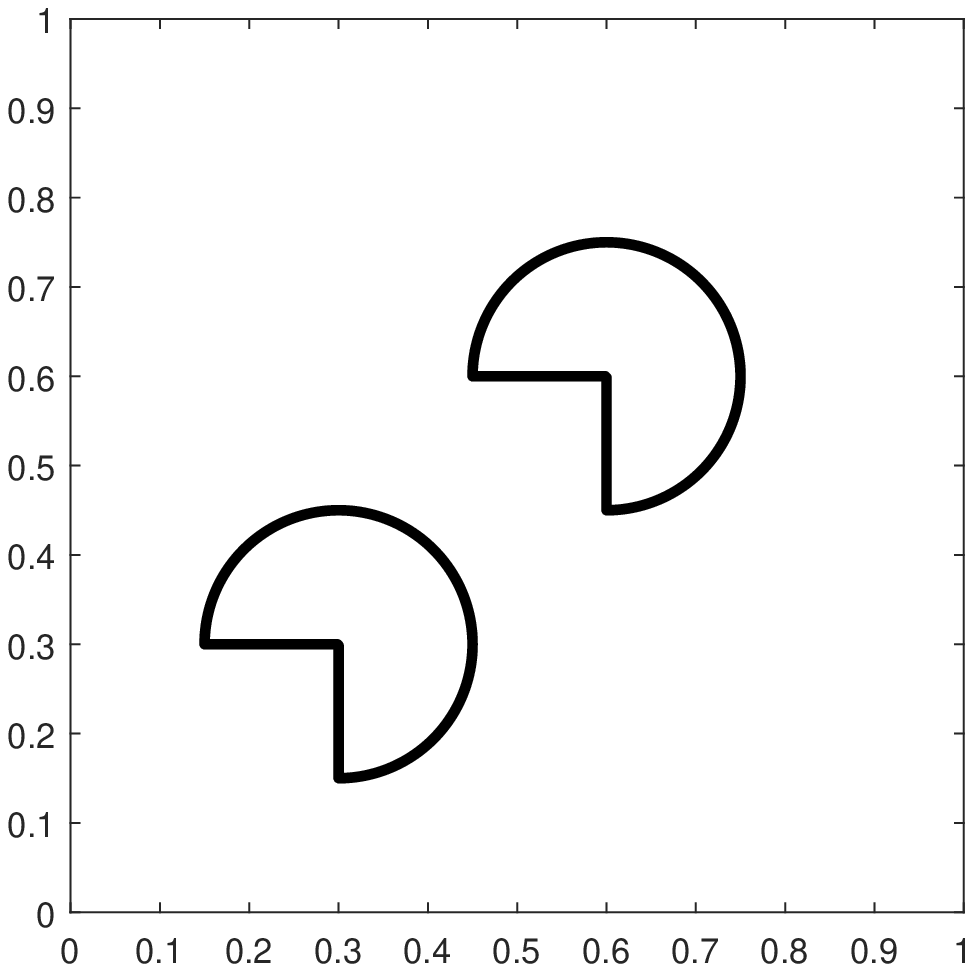}}
	\subfloat[$E_z$]{\includegraphics[scale=0.4]{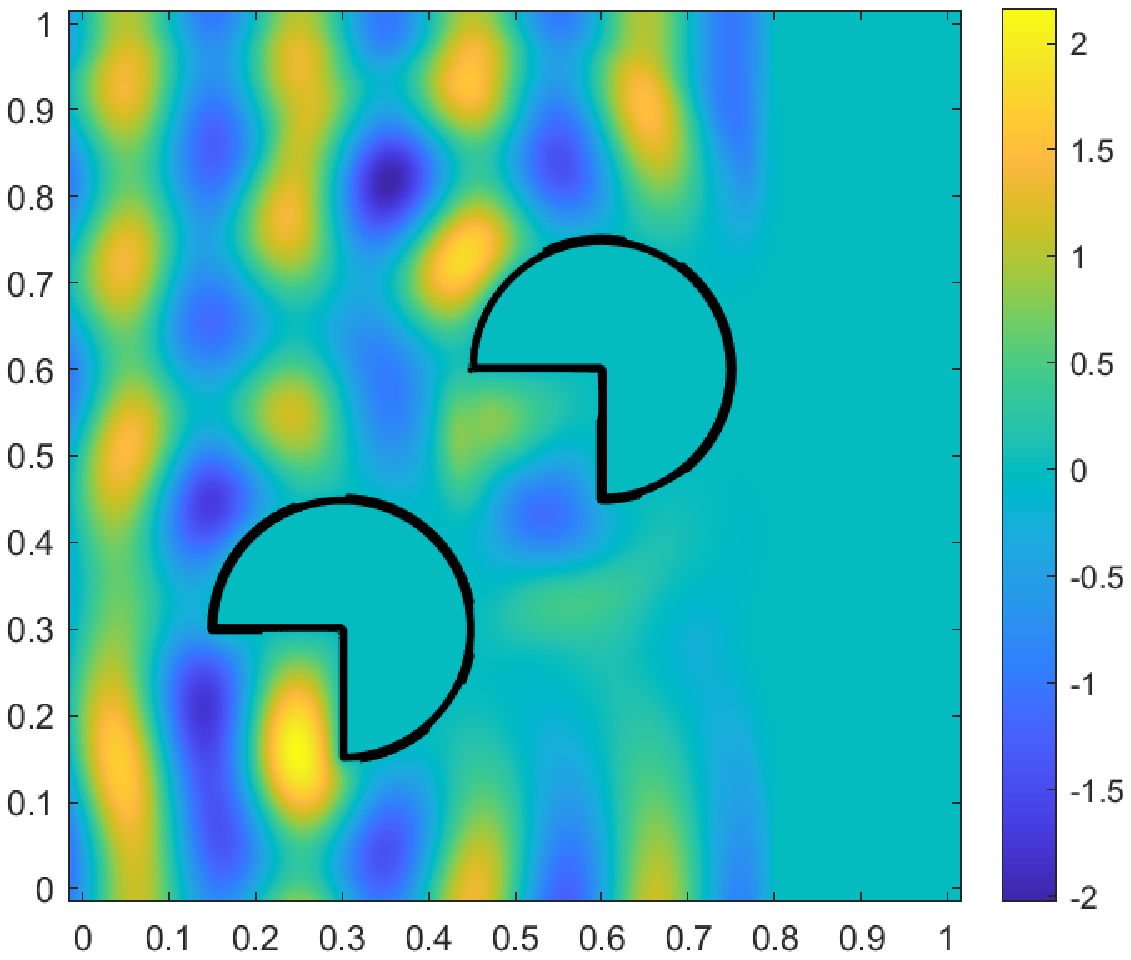}}		
	\caption{Solution snapshot of $E_z$ at time $T = 0.8$ with the sine incident wave and two circular wedge-shaped PECs}
	\label{fig:corner_double_pizza_PEC}
\end{figure}}

	\subsection{ Effect of CFL numbers on long time dynamic of  $E_z$} 
\begin{figure}[H]
	\centering
	\includegraphics[scale=0.2]{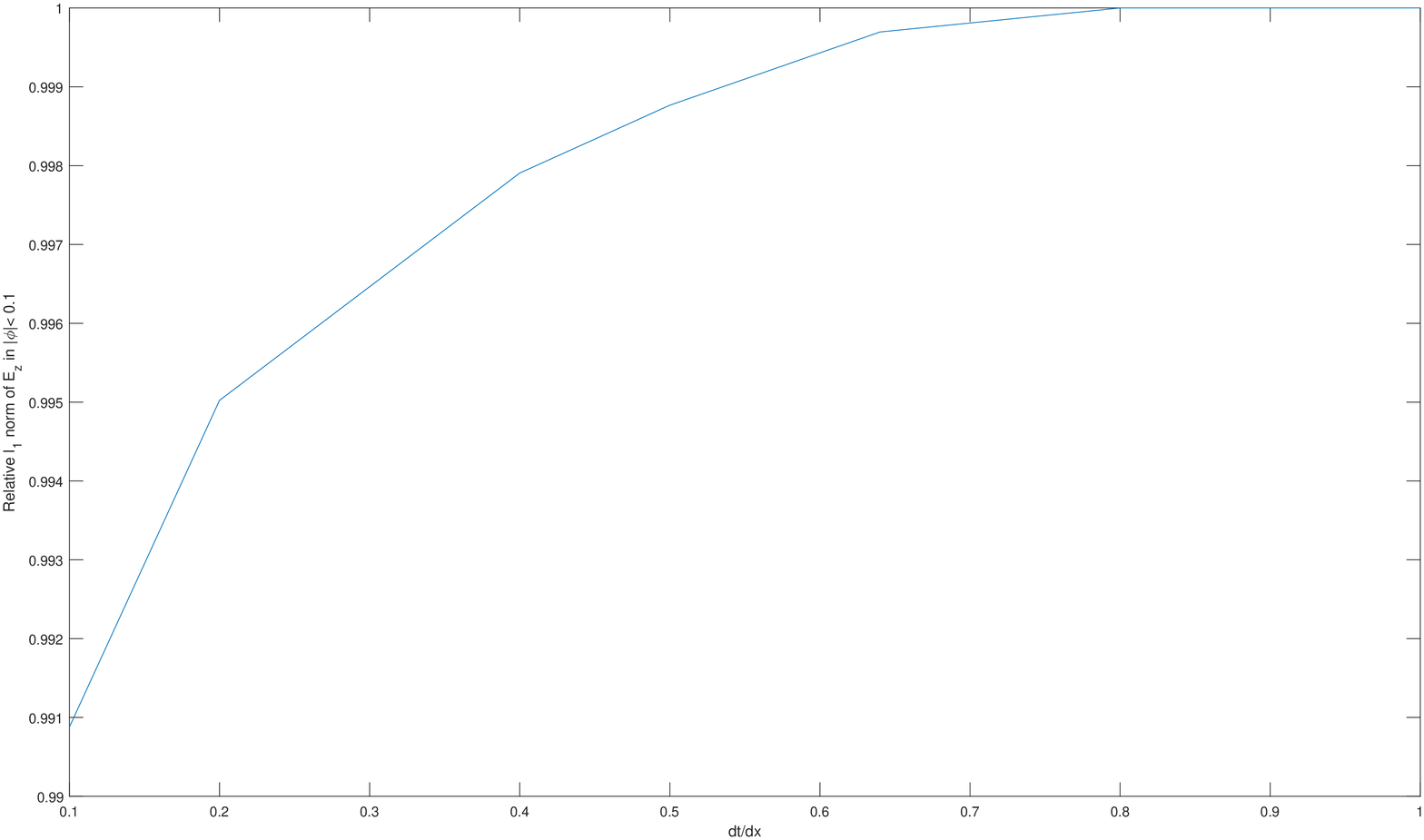}
	\caption{Relative $l_1$ norm of $E_z$ with respect to different CFL numbers $\Delta t/\Delta x$}
	\label{fig:l1_cfl_vary}
\end{figure}

\begin{figure}[H]
		\centering
		\includegraphics[scale=0.2]{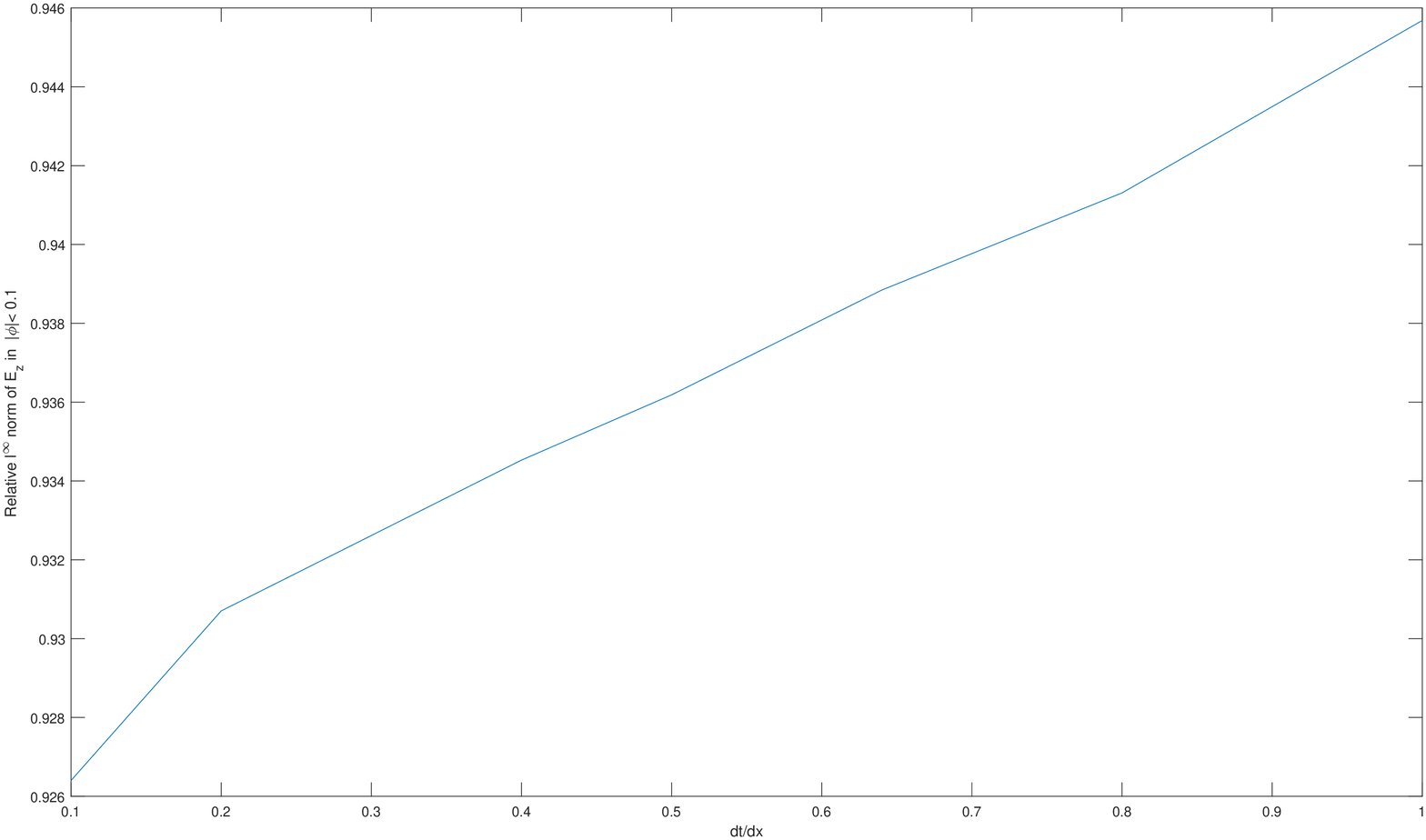}
		\caption{Relative $l^\infty$ norm of $E_z$ with respect to different CFL numbers $\Delta t/ \Delta x$}
		\label{fig:linfty_cfl_vary}
	\end{figure}

Following the discussion in \cite{zl21} we investigate how the dissipative behaviors of our method depend on the CFL numbers $\Delta t / \Delta x$. First we analyze what happens to the peaks of $E_z$ as the CFL numbers vary. Returning to the setup of our numerical experiment for the circular PEC with the plane incident wave, we compute the numerical solutions on the \emph{fixed} spatial mesh size of $\Delta x = \Delta y = \frac{1}{160}$ while the CFL numbers are chosen to be $0.1, 0.2, 0.4, 0.64, 0.8, 1$. The computations are carried out until the terminal time $T = 0.8$. The solutions obtained are then compared with the ``reference" solution computed on the mesh size of $\Delta x = \Delta y = \frac{1}{640}$ with the CFL number equal to $1$. The comparisons are made based on the $l^\infty$ and $l_1$ norms of the solutions over the collar of the width $0.1$ outside the PEC. The results are displayed in Figures \ref{fig:l1_cfl_vary} and \ref{fig:linfty_cfl_vary}. 

	We see that the amplitudes remain sharp, retaining close to $99\%$ and $92\%$ of the averaged and the maximum amplitudes, respectively, for all CFL numbers tested. We note though a larger decay of the peaks as the CFL numbers decrease.

\begin{figure}[H]
	\centering
	\includegraphics[scale=0.2]{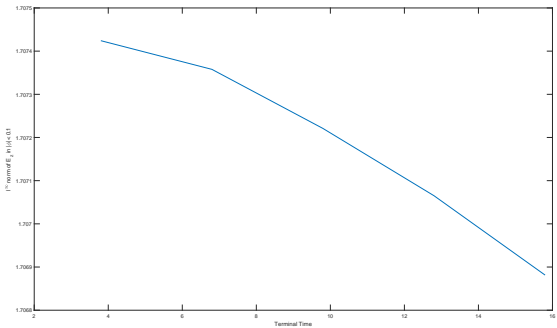}
	\caption{$l^\infty$ norm of $E_z$ with respect to different terminal times}
	\label{fig:longtime}
\end{figure}

We finally examine whether the peaks of $E_z$ are preserved in  the long time simulations. We choose the spatial mesh of $\Delta x = \Delta y = \frac{1}{160}$ with the CFL number set to $1$. We again use the circular PEC with the plane incident wave, and compute the solutions until the terminal time $T = 3.8,6.8,9.8, 12.8$. These correspond to spanning the distance of roughly 40 times the wavelength of the incident wave. Using the $l^\infty$ norm of $E_z$ as the measurement we see from Figure \ref{fig:longtime} that the amplitudes decay less than $0.1\%$ as the terminal time advances from $T = 3.8$ to $T= 12.8$.    

\section{Conclusion}
\label{sec:conc}

Our finite difference method produces second order accurate solutions on uniform orthogonal grids in which general PECs are embedded with their level set representations. The method is straightforward to implement relying on repeated applications of first order finite difference schemes based on 5-point stencils. It is a simplified, yet as effective modification of the existing algorithm on the point shifted grid \cite{zl21}. In our future work we wish to apply our method to three dimensional problems as well as moving PEC objects.

\section*{Acknowledgement}
HL is indebted to Haiyu Zou for providing and explaining the codes which were used to generate the numerical results in https://arxiv.org/abs/\\ 2104.14675.

\bibliographystyle{siamplain}
\bibliography{references}
\end{document}

%% file: ex_shared.tex
\usepackage{lipsum}
\usepackage{amsfonts}
\usepackage{graphicx}
\usepackage{epstopdf}
\usepackage{algorithmic}
\ifpdf
  \DeclareGraphicsExtensions{.eps,.pdf,.png,.jpg}
\else
  \DeclareGraphicsExtensions{.eps}
\fi

\newsiamremark{remark}{Remark}
\newsiamremark{hypothesis}{Hypothesis}
\crefname{hypothesis}{Hypothesis}{Hypotheses}
\newsiamthm{claim}{Claim}

\headers{A ghost-point based second order finite difference method}{H. Lee and Y. Liu}

\title{A ghost-point based second order accurate finite difference method on uniform orthogonal grids for electromagnetic scattering around PEC\thanks{Submitted to the editors DATE.}}

\author{Hwi Lee\thanks{School of Mathematics, Georgia Institute of Technology, Atlanta, GA 30332 USA.
  (\email{hlee995@gatech.edu}, \email{yingjie@math.gatech.edu}).}
\and Yingjie Liu \footnotemark[2]}

\usepackage{amsopn}


%% file: BFECC_uniform_grid_manuscript.bbl
\begin{thebibliography}{10}

\bibitem{AS98}
{\sc D.~Adalsteinsson and J.~A. Sethian}, {\em The {F}ast {C}onstruction of
  {E}xtension {V}elocities in {L}evel {S}et {M}ethods}, J. Comput. Phys., 148
  (1999), pp.~2--22.

\bibitem{a18}
{\sc J.~B. Angel, J.~W. Banks, and W.~D. Henshaw}, {\em High-order upwind
  schemes for the wave equation on overlapping grids: {M}axwell's equations in
  second-order form}, Journal of Computational Physics, 352 (2018),
  pp.~534--567.

\bibitem{ch01}
{\sc D.~L. Chopp}, {\em Some improvements of the fast marching method}, SIAM
  Journal on Scientific Computing, 23 (2001), pp.~230--244.

\bibitem{d13}
{\sc M.~Detrixhe, F.~Gibou, and C.~Min}, {\em A parallel fast sweeping method
  for the {E}ikonal equation}, Journal of Computational Physics, 237 (2013),
  pp.~46--55.

\bibitem{DDH01}
{\sc A.~Ditkowski, K.~Dridi, and J.~S. Hesthaven}, {\em Convergent {C}artesian
  grid methods for {M}axwell’s equations in complex geometries}, J. Comput.
  Phys., 170 (2001), pp.~39--80.

\bibitem{dupontBackForthError2003}
{\sc T.~F. Dupont and Y.~Liu}, {\em Back and forth error compensation and
  correction methods for removing errors induced by uneven gradients of the
  level set function}, Journal of Computational Physics, 190 (2003),
  pp.~311--324.

\bibitem{dupontliu07}
{\sc T.~F. Dupont and Y.~Liu}, {\em Back and forth error compensation and
  correction methods for semi-lagrangian schemes with application to level set
  interface computations}, Mathematics of Computation,  (2007), pp.~647--668.

\bibitem{fangLocallyConformedFinitedifference1993}
{\sc J.~Fang and J.~Ren}, {\em A locally conformed finite-difference
  time-domain algorithm of modeling arbitrary shape planar metal strips}, IEEE
  transactions on microwave theory and techniques, 41 (1993), pp.~830--838.

\bibitem{Fedkiw99}
{\sc R.~P. Fedkiw, T.~Aslam, B.~Merriman, and S.~Osher}, {\em A non-oscillatory
  {E}ulerian approach to interfaces in multimaterial flows (the ghost fluid
  method)}, J. Comput. Phys., 152 (1999), pp.~457--492.

\bibitem{f02}
{\sc N.~Flyer and P.~N. Swarztrauber}, {\em The convergence of spectral and
  finite difference methods for initial-boundary value problems}, SIAM Journal
  on Scientific Computing, 23 (2002), pp.~1731--1751.

\bibitem{h16}
{\sc L.~Hu, Y.~Li, and Y.~Liu}, {\em A limiting strategy for the back and forth
  error compensation and correction method for solving advection equations},
  Mathematics of Computation, 85 (2016), pp.~1263--1280.

\bibitem{jackson1999classical}
{\sc J.~D. Jackson}, {\em Classical electrodynamics}, 1999.

\bibitem{k05}
{\sc B.~Kim, Y.~Liu, I.~Llamas, and J.~R. Rossignac}, {\em {F}lowfixer: Using
  {BFECC} for fluid simulation}, tech. report, Georgia Institute of Technology,
  2005.

\bibitem{k07}
{\sc D.~Komatitsch and R.~Martin}, {\em An unsplit convolutional perfectly
  matched layer improved at grazing incidence for the seismic wave equation},
  Geophysics, 72 (2007), pp.~SM155--SM167.

\bibitem{kpy04}
{\sc H.-O. Kreiss, N.~A. Petersson, and J.~Ystr{\"o}m}, {\em Difference
  approximations of the {N}eumann problem for the second order wave equation},
  SIAM Journal on Numerical Analysis, 42 (2004), pp.~1292--1323.

\bibitem{k93}
{\sc K.~S. Kunz and R.~J. Luebbers}, {\em The finite difference time domain
  method for electromagnetics}, CRC press, 1993.

\bibitem{ln21}
{\sc Y.-M. Law and J.-C. Nave}, {\em {FDTD} schemes for {Maxwell’s}
  {E}quations with {E}mbedded {P}erfect {E}lectric {C}onductors {B}ased on the
  {C}orrection {F}unction {M}ethod}, Journal of Scientific Computing, 88
  (2021), pp.~1--28.

\bibitem{liuOverlappingYeeFDTD2009}
{\sc J.~Liu, M.~Brio, and J.~V. Moloney}, {\em Overlapping {Yee} {FDTD}
  {Method} on {Nonorthogonal} {Grids}}, Journal of Scientific Computing, 39
  (2009), pp.~129--143.

\bibitem{m10}
{\sc C.~Min}, {\em On reinitializing level set functions}, Journal of
  computational physics, 229 (2010), pp.~2764--2772.

\bibitem{m01}
{\sc M.~Min and C.~Teng}, {\em The instability of the {Y}ee scheme for the
  “magic time step”}, Journal of Computational Physics, 166 (2001),
  pp.~418--424.

\bibitem{m97}
{\sc J.~P. Morris, P.~J. Fox, and Y.~Zhu}, {\em Modeling low {R}eynolds number
  incompressible flows using {SPH}}, Journal of computational physics, 136
  (1997), pp.~214--226.

\bibitem{m81}
{\sc G.~Mur}, {\em The modeling of singularities in the finite-difference
  approximation of the time-domain electromagnetic-field equations}, IEEE
  Transactions on Microwave Theory and Techniques, 29 (1981), pp.~1073--1077.

\bibitem{o01}
{\sc S.~Osher and R.~P. Fedkiw}, {\em Level set methods: an overview and some
  recent results}, Journal of Computational physics, 169 (2001), pp.~463--502.

\bibitem{o88}
{\sc S.~Osher and J.~A. Sethian}, {\em Fronts propagating with
  curvature-dependent speed: {A}lgorithms based on {H}amilton-{J}acobi
  formulations}, Journal of computational physics, 79 (1988), pp.~12--49.

\bibitem{jin-fa-leeModelingThreedimensionalDiscontinuities1992}
{\sc R.~Palandech, R.~Mittra, et~al.}, {\em Modeling three-dimensional
  discontinuities in waveguides using nonorthogonal {FDTD} algorithm}, IEEE
  Transactions on Microwave Theory and Techniques, 40 (1992), pp.~346--352.

\bibitem{Peng99}
{\sc D.~Peng, B.~Merriman, S.~Osher, H.~Zhao, and M.~Kang}, {\em A {PDE}-based
  fast {L}ocal {L}evel {S}et {M}ethod}, J. Comput. Phys., 155 (1999),
  pp.~410--438.

\bibitem{r99}
{\sc C.~J. Railton and J.~B. Schneider}, {\em An analytical and numerical
  analysis of several locally conformal {FDTD} schemes}, IEEE Transactions on
  Microwave Theory and Techniques, 47 (1999), pp.~56--66.

\bibitem{SFKLR08}
{\sc A.~Selle, R.~Fedkiw, B.~Kim, Y.~Liu, and J.~Rossignac}, {\em An
  {U}nconditionally {S}table {M}accormack {M}ethod}, J. Sci. Comput., 35
  (2008), pp.~350--371.

\bibitem{sussmanLevelSetApproach1994}
{\sc M.~Sussman, P.~Smereka, and S.~Osher}, {\em A level set approach for
  computing solutions to incompressible two-phase flow}, Journal of
  Computational physics, 114 (1994), pp.~146--159.

\bibitem{ta05}
{\sc A.~Taflove, S.~C. Hagness, and M.~Piket-May}, {\em Computational
  electromagnetics: the finite-difference time-domain method}, The Electrical
  Engineering Handbook, 3 (2005), pp.~629--670.

\bibitem{umashankarCalculationExperimentalValidation1987a}
{\sc K.~Umashankar, A.~Taflove, and B.~Beker}, {\em Calculation and
  experimental validation of induced currents on coupled wires in an arbitrary
  shaped cavity}, IEEE Transactions on Antennas and Propagation, 35 (1987),
  pp.~1248--1257.

\bibitem{wangBackForthError2019}
{\sc X.~Wang and Y.~Liu}, {\em Back and forth error compensation and correction
  method for linear hyperbolic systems with application to the {Maxwell}'s
  equations}, Journal of Computational Physics: X, 1 (2019), p.~100014.

\bibitem{yeeNumericalSolutionInital1966}
{\sc K.~Yee}, {\em Numerical solution of inital boundary value problems
  involving {M}axwell's equations in isotropic media}, IEEE Transactions on
  antennas and propagation, 14 (1966), pp.~302--307.

\bibitem{y92}
{\sc K.~S. Yee, J.~S. Chen, and A.~H. Chang}, {\em Conformal finite difference
  time domain (fdtd) with overlapping grids}, in IEEE Antennas and Propagation
  Society International Symposium 1992 Digest, IEEE, 1992, pp.~1949--1952.

\bibitem{zhao05}
{\sc H.~Zhao}, {\em A fast sweeping method for {E}ikonal equations},
  Mathematics of computation, 74 (2005), pp.~603--627.

\bibitem{zl21}
{\sc H.~Zou and Y.~Liu}, {\em A finite difference method on irregular grids
  with local second order ghost point extension for solving {M}axwell's
  equations around curved {PEC} objects}, Journal of Computational Physics, 463
  (2022), p.~111273.

\end{thebibliography}
